\documentclass[12pt]{article}

\usepackage{amsthm}

\usepackage[left=1cm,right=1cm,top=0.5cm,bottom=1.5cm,bindingoffset=0cm]{geometry}

\usepackage[T2A]{fontenc}
\usepackage[utf8]{inputenc}
\usepackage[english]{babel}
\usepackage{amsmath}
\usepackage{amsfonts}
\usepackage{mathtools}
\usepackage{graphicx}
\usepackage{systeme}
\usepackage{tikz}
\usepackage{blindtext}
\usepackage{amssymb}
\usepackage{hyperref}
\usepackage{csquotes}
\usepackage{bbm}
\usepackage[
backend=biber,
style=numeric,
sorting=none
]{biblatex}
\usepackage[normalem]{ulem}
\graphicspath{ {./} }

\newtheorem{theorem}{Theorem}

\newtheorem{lemma}{Lemma}[section]
\newtheorem*{lemma*}{Lemma}
\newtheorem*{den*}{Denotion}

\newtheorem*{remark*}{Remark}
\newtheorem{construction}{Construction}

\theoremstyle{remark}

\theoremstyle{definition}
\newtheorem{definition}{Definition}[section]

\newcommand{\mP}{\mathbf{P}}

\title{Stability of the independence number of $G(n, r, 1)$ graphs\footnote{The research is supported by the grant of Theoretical Physics and Mathematics Advancement Foundation “BASIS”.}}
\date{}
\author{M.M. Koshelev\footnote{Moscow Institute of Physics and Technology, Moscow State University, mkoshelev99@gmail.com}}

\begin{document}

\maketitle

\sloppy

\begin{abstract}
    In this paper we obtain the stability theorem for the independence number of $G(n, r, 1)$ graphs. This result was previously stated in the paper of M. Pyaderkin but the proof there was incorrect. We introduce the correct proof of the key lemma and thus finally complete the proof of this theorem.
\end{abstract}

\section{Introduction}

This note is devoted to the study of distance graphs $G(n,r,s)$, also known as Johnson graphs. Let us give a definition of these graphs.
$$G(n, r, s) :=(V, E),~V=\binom{[n]}{r}, ~ E=\{(A, B):~ |A\cap B| =s\},$$
that is, the vertices of the graph are all possible $r$-element subsets of $[n] := \{1, \ldots, n\}$, and edges connect pairs of sets that intersect in exactly $s$ elements. In addition, we will use the notation $G_p(n,r,s)$ for a random graph obtained from $G(n,r,s)$ by deleting each of its edges independently of the others with probability $1 - p$.

A vast amount of modern scientific literature in the field of extremal combinatorics is devoted to $G(n,r,s)$ and $G_p(n,r,s)$ graphs (see, for example, \cite{Codes, RaiMeta, Rai1, MSp, RSFW, KRdGn}). One of the mainstream directions here is the study of independence (see \cite{EKR, FrW, FrFu}) and chromatic (see \cite{Hyper, ZChrom}) numbers of Johnson graphs. One might also refer to the survey \cite{RaiMetaAlphaChrom} devoted to this field of study.

Recall the definition of the independence number.

\begin{definition}
    The independence number of a graph $G$ is the maximum number of vertices of $G$ such that the subgraph induced on these vertices is empty. The standard denotation for the independence number is $\alpha(G)$.
\end{definition}

Speaking of the study of independence numbers of Johnson graphs, one cannot fail to mention the series of works \cite{BNR, BBN, Das}, which culminated in the work \cite{BKL}, where a theorem on the stability of the independence numbers of graphs $G(n,r,0)$, also known as Kneser graphs, was proved.

\begin{theorem}
    Let $\varepsilon > 0$ be an arbitrary number, and let $n \geq 2r+1$. Then for every $p \geq (1 + \varepsilon)p_0$, where 
    $$
    p_0 = \begin{cases}
        3/4, & n = 2r + 1, \\
        \log\left(n\binom{n-1}{r}\right) / \binom{n-r-1}{r-1},& n > 2r+1,
    \end{cases}
    $$
    the graph $G_p(n,r,0)$ will with high probability (w.h.p.) have independence number equal to $\binom{n-1}{r-1}$.
    
    Moreover the equality $\alpha(G_p(n,r,0)) = \binom{n-1}{r-1}$ w.h.p. does not hold for any $p \leq (1 - \varepsilon)p_0$.  
    
    We will use the notation w.h.p. hereinafter for a sequence of events $A_n$ such that $\mP(A_n) \to 1$.
\end{theorem}

Let us discuss this result. First of all, it is not obvious that $\alpha(G(n,r,0)) = \binom{n-1}{r-1}$, let alone $\alpha(G_p(n,r,s))$. This is the famous Erd\H os--Ko--Rado theorem (see \cite{EKR}), a beautiful example of probabilistic method in combinatorics. We will not go into details of the proof of the hard part of the theorem, the $\alpha(G(n,r,0)) \leq \binom{n-1}{r-1}$ inequality, but will give an example of an independent set of size $\binom{n-1}{r-1}$. The example of such a set (and in fact the only possible example) is the subgraph that we will call a star with center $x$, the set of vertices $S_x = \{v\in G(n,r,0)|x \in v\}$. Note that the name star is sometimes used in graph theory to denote a tree of a certain type. The stars indeed play a crucial role in the study of $\alpha(G(n,r,s))$. In fact, the celebrated theorem of Frankl and F\"uredi (see \cite{FrFu}) states that as long as $r \geq 2s + 1$ and $n$ is big enough the independence number of $G(n,r,s)$ is equal to $\binom{n-s-1}{r-s-1}$ and is attained at a star $S_{x_1,\dotsc, x_{s+1}} := \{v \in G(n,r,s)| \{x_1,\dotsc, x_{s+1}\} \subseteq v\}$.

Of course, it is natural to ask a similar question about the stability of the independence number for graphs $G(n,r,s)$ with $s > 0$. Unfortunately, very little is known about this case. In the works \cite{Pyad1Eng,Pyad2Eng}, the stability of the independence number in order of magnitude was studied, and the only work known to the author on the exact stability of the independence number of graphs $G(n,r,s)$ with $s > 0$ is the work \cite{PyadBug}. In this paper the following result is claimed.

\begin{theorem}
    \label{thGnr1}
    Let $r \geq 4$. Then w.h.p. $\alpha(G_{1/2}(n,r,1)) = \binom{n-2}{r-2}$.
\end{theorem}

One of the key lemmas in the proof of Theorem \ref{thGnr1} is the following lemma.

\begin{lemma}
\label{bugged}
    Let $r \geq 4$. There exists some $t_0 \in (1/2, 1)$ such that w.h.p. any $A$ satisfying the conditions
     \begin{itemize}
         \item  $|A| > \binom{n-2}{r-2}$;

         \item $d(A) > \binom{t_0n-2}{r-2}$ (here and henceforth, by $d(A)$ we mean $\max_{1 \leq i < j \leq n} |A\cap S_{i,j}|$);
     \end{itemize}
    is not independent in $G_{1/2}(n,r,1)$.
\end{lemma}

Unfortunately, the proof of Lemma \ref{bugged} from the paper \cite{PyadBug} is incorrect. In fact, the paper proves the following:

\begin{lemma}
\label{bug}
    Let $r \geq 4$. There exists some $t_0 \in (1/2, 1)$ such that w.h.p. any $A$ satisfying the conditions
     \begin{itemize}
         \item  $|A| > \binom{n-2}{r-2}$;

         \item $d(A) > \binom{t_0n-2}{r-2}$;

         \item There exists a vertex $u \in A: |\{i,j\} \cap u| = 1$ (where $i,j$ is the center of the star having the maximum intersection with $A$)
     \end{itemize}
    is not independent in $G_{1/2}(n,r,1)$.
\end{lemma}

Thus, to correct the proof, it is sufficient to prove the following lemma.

\begin{lemma}
\label{fix}
     Let $r \geq 4$. There exists some $t_0 \in (1/2, 1)$ such that w.h.p. any $A$ satisfying the conditions
     \begin{itemize}
         \item  $|A| = \binom{n-2}{r-2} + 1$;

         \item $d(A) > \binom{t_0n-2}{r-2}$;

         \item all vertices $u \in A$ either contain both $i$ and $j$, or are subsets of the set $[n] \setminus \{i, j\}$ (where $i,j$ is the center of the star having the maximum intersection with $A$)
     \end{itemize}
    is not independent in $G_{1/2}(n,r,1)$.
\end{lemma}

In this paper we present the proof of Lemma \ref{fix} thus completing the proof of Theorem \ref{thGnr1}.

\section{Proof of Lemma \ref{fix}}

Before proceeding to the proof of Lemma \ref{fix}, let us introduce some auxiliary notations and conventions. Everywhere, unless explicitly stated otherwise, we will assume that for a set $A$ a star with maximum intersection with $A$ is the star $S_{n-1,n}$. Let $X = X(A)$ be the set of vertices of $A$ not belonging to $S_{n-1,n}$, $|X| = x$. Furthermore, let $I(X)$ be the set of numbers from $[n-2]$ that are contained in at least one vertex of the set $X$, $|I(X)| = i(X)$. Finally, throughout this text, we assume that $\binom{a}{b} = 0$ whenever $a < b$ or $b < 0$.

For a proof of one of the cases of Lemma \ref{fix}, we will need an estimate for the number of edges in subgraphs of $G(n,r,1)$. It follows trivially from Theorems 4 and 5 in \cite{ShKup} (see also Theorem 7 from \cite{Shu1}). We formulate this statement in a form convenient for us.

\begin{theorem}
    \label{cntE}
    Let $r \geq 3$, and $l > (1 + \varepsilon)\binom{n}{r-2}$ for some $\varepsilon > 0$. Then for any set of vertices $L, |L| = l$ of the graph $G(n,r,1)$, the following relation holds.
    $$
        e(L) = \Theta\left(\frac{l^2}{n}\right),
    $$
    where $e(L)$ is the number of edges within $L$.
\end{theorem}

In addition to this, we will need the following simple estimate for the difference of binomial coefficients.

\begin{lemma}
    \label{Tech} 
    Let $n, r, i$ be arbitrary integers, $i \leq n, n \geq 2r - 4$. Then there exists a $c = c(r) > 0$ such that
    $$
        \binom{n-r-1}{r-2}-\binom{n-r-1-i}{r-2} \geq cin^{r-3}.
    $$
\end{lemma}

Next, we will need the following well-known estimate due to Chernoff (see \cite{Chernoff})

\begin{theorem}
    \label{Chernoff} 
    Let $X \sim Bin(n, p), \mu = np$. Then 
    $$
        \mP(X > (1 + \delta)\mu) \leq e^{-\frac{\delta^2}{2 + \delta}\mu}.
    $$
\end{theorem}

Finally, in the proof we will use the following construction.

\begin{construction}
    Consider some set $A$ of the vertices of $G(n,r,1)$. Let $I(X) = I(X(A)) = \{l_1, \dotsc, l_{i(X)}\}$. Consider vertices $u_1, \dotsc, u_{i(X)}$ such that $u_{j} \in X, l_{j} \in u_j$. We define the sets $B_j, j = 1, \dotsc, i(X)$ as sets of vertices with the following two conditions:
        \begin{itemize}
            \item For any $v \in B_j$,  $\{l_j, n-1, n\} \subset v$.

            \item For any $v \in B_j$,  $(\{l_1, \dotsc, l_{j-1}\}\cup u_j \setminus \{l_j\}) \cap v = \varnothing$.
        \end{itemize}
        It is easy to see that $$|B_j| \geq \binom{n-r-j-1}{r-3}$$ and, furthermore, all $B_j$ are disjoint. Thus it follows that $$|B_1 \sqcup \dotsc \sqcup B_{i(X)}| \geq \sum_{j=1}^{i(X)} \binom{n-r-j-1}{r-3}.$$ Besides, it is easy to see that all vertices from $B_i$ intersect with $u_i$ by exactly one element (and, consequently, there is an edge between them).
\end{construction}

We are ready to proceed with the proof. We will split it into two lemmas, from the combination of which Lemma \ref{fix} follows trivially.

\begin{lemma}
    \label{fixBipartite}
    Let $r \geq 4$. Consider all possible sets $A$ with the following conditions:
    \begin{enumerate}
        \item  $|A| = \binom{n-2}{r-2} + 1$;

         \item $d(A) > \binom{t_0n-2}{r-2}$;

         \item all vertices $u \in A$ either contain both $i$ and $j$, or are subsets of the set $[n] \setminus \{i, j\}$ (where $i,j$ is the center of the star having the maximum intersection with $A$);

         \item $\sum_{j = 1}^{i(X)} \binom{n-j-1-r}{r-3} \geq 3x$. 
    \end{enumerate}

    For any $t_0 \in (1/2, 1)$ w.h.p. no such set is independent in $G_{1/2}(n,r,1)$.
\end{lemma}

\begin{lemma}
\label{fixTuran}
    Let $r \geq 4$. Consider all possible sets $A$ with the following conditions:
    \begin{enumerate}
        \item  $|A| = \binom{n-2}{r-2} + 1$;

         \item $d(A) > \binom{t_0n-2}{r-2}$;

         \item all vertices $u \in A$ either contain both $i$ and $j$, or are subsets of the set $[n] \setminus \{i, j\}$ (where $i,j$ is the center of the star having the maximum intersection with $A$);

         \item $\sum_{j = 1}^{i(X)} \binom{n-j-1-r}{r-3} < 3x$.
    \end{enumerate}

    There exists $t_0 \in (1/2, 1)$ such that w.h.p. no such set is independent in $G_{1/2}(n,r,1)$.
\end{lemma}

The rest of the paper is organized as follows. In section \ref{TechPr}, we prove Lemma \ref{Tech}. Section \ref{BipPr} contains the proof of Lemma \ref{fixBipartite}, and section \ref{TurPr} contains the proof of Lemma \ref{fixTuran}.

\subsection{Proof of Lemma \ref{Tech}}
\label{TechPr}

First assume $i < n / 2$. We use the following equality:
$$
    \binom{n-r-1}{r-2} = \sum_{j=0}^{r-2}\binom{i}{j}\binom{n-r-1-i}{r-2-j} \geq \binom{n-r-1-i}{r-2}+i\binom{n-r-1-i}{r-3},
$$
from which
$$
    \binom{n-r-1}{r-2}-\binom{n-r-1-i}{r-2} \geq i\binom{n-r-i-1}{r-3} \geq c_1in^{r-3}
$$
for some $c_1 > 0$ when $i \leq n/2$.

For $i \geq n/2$ we have $$\binom{n-r-1-i}{r-2} \leq (1-c_2)\binom{n-r-1}{r-2},$$ from which the corresponding difference is not less than $$c_2\binom{n-r-1}{r-2} > c_3n^{r-2} \geq c_3in^{r-3}$$ for some $c_3 > 0$.

\subsection{Proof of Lemma \ref{fixBipartite}}
\label{BipPr}

Consider an arbitrary $A$ satisfying the conditions of Lemma \ref{fixBipartite} and construct from it a subset of vertices $Ess(A)$ as follows: first, we put all vertices from $A \setminus X(A)$ in $Ess(A)$. Furthermore, if $I(X) = \{l_1, \dotsc, l_{i(X)}\}$, then we place in $Ess(A)$ vertices $v_1, \dotsc, v_{i(X)}$ with the conditions $l_k \in v_{k} \in A$ (some of the $v_i$ may, of course, coincide). Of course, $Ess(A)$ is not uniquely determined by $A$, but we do not need this. It is only important for us that w.h.p. none of the possible sets $Ess(A)$ are independent, and, therefore, the same is true for all possible $A$. It remains to prove that w.h.p. none of the $Ess(A)$ are independent in $G_{1/2}(n,r,1)$. To do this, we first fix the star and $i(X) \leq n$ vertices from $Ess(A)$ outside the star. As we know, in the star there are at least $$\sum_{j = 1}^{i(X)} \binom{n-j-1-r}{r-3}$$ vertices that have edges to the fixed vertices. Of these, we must include all but at most $x$ vertices in the set $Ess(A)$. That is, in order to be able to form at least one independent $Ess(A)$ with the chosen $i(X)$ vertices and the chosen star, at least $$\frac{2}{3}\sum_{j = 1}^{i(X)} \binom{n-j-1-r}{r-3}$$ vertices of the star must have their edges to the chosen vertices disappear. But for each specific vertex, these edges will disappear with a probability of at most $1/2$. This means that out of $$\sum_{j = 1}^{i(X)} \binom{n-j-1-r}{r-3}$$ variables with a $Bern(p), p \leq 1/2$ distribution, at least $$\frac{2}{3}\sum_{j = 1}^{i(X)} \binom{n-j-1-r}{r-3}$$ must be equal to 1. By Theorem \ref{Chernoff}, the probability of such an event does not exceed 
$$
    e^{-\frac{\frac{1}{6^2}}{2+\frac{1}{6}}\frac{1}{2}\sum_{j = 1}^{i(X)} \binom{n-j-1-r}{r-3}} = e^{-\frac{1}{156}\sum_{j = 1}^{i(X)} \binom{n-j-1-r}{r-3}} = e^{-\frac{1}{156}\left(\binom{n-r-1}{r-2}-\binom{n-r-1-i(X)}{r-2}\right)}.
$$
Then the probability of the existence of at least one independent $Ess(A)$ is no more than
$$
    \sum_{i=1}^{n-2}\binom{n}{2}\binom{n}{r}^{i}e^{-\frac{1}{156}\left(\binom{n-r-1}{r-2}-\binom{n-r-1-i}{r-2}\right)} \leq \sum_{i=1}^{n-2}e^{2\ln n + ir\ln n -\frac{1}{156}\left(\binom{n-r-1}{r-2}-\binom{n-r-1-i}{r-2}\right)}.
$$
Using Lemma \ref{Tech}, we finally obtain
$$
    \sum_{i=1}^{n-2}e^{2\ln n + ir\ln n -\frac{1}{156}\left(\binom{n-r-1}{r-2}-\binom{n-r-1-i}{r-2}\right)} \leq \sum_{i=1}^{n-2}e^{2\ln n - c_5in^{r-3}} = o(1).
$$

\subsection{Proof of Lemma \ref{fixTuran}}
\label{TurPr}

Consider the vertices from $X(A)$. These are vertices of the graph $G(i(X), r, 1)$. Let us prove that the conditions of the lemma imply the inequality $$x > 2\binom{i(X)-2}{r-2}.$$ Indeed, it is easy to see that under the conditions of the lemma
$$
    x > \frac{1}{3}\sum_{j = 1}^{i(X)} \binom{n-j-1-r}{r-3} = \frac{1}{3}\left(\binom{n-r-1}{r-2}-\binom{n-r-1-i(X)}{r-2}\right).
$$
From this, by virtue of Lemma \ref{Tech}, the inequality $x > ci(X)n^{r-3}$ immediately follows. At the same time $$\binom{i(X)-2}{r-2} < i^{r-2}(X).$$ Furthermore, we note since of $$x \leq \binom{n-2}{r-2}-\binom{t_0n-2}{r-2},$$ we can choose $t_0$ such that the inequality $x < \varepsilon n^{r-2}$ holds for any given $\varepsilon > 0$. From this, the inequality $$ci(X)n^{r-3} < \varepsilon n^{r-2} \Leftrightarrow i(X) < \varepsilon'n$$ will follow for any given $\varepsilon' > 0$. Finally, from this inequality it follows that $$2i^{r-2}(X) < 2(\varepsilon')^{r-3}n^{r-3}i(X).$$ It remains to take $$\varepsilon' = \sqrt[r-3]{\frac{c}{2}}$$ to obtain the required statement.

In addition, it is easy to see that by virtue of the inequalities $$ci(X)n^{r-3} < x \leq \binom{i(X)}{r},$$ the relation $$i(X) > \Tilde{c}n^{\frac{r-3}{r-1}}$$ holds.

From this it trivially follows that $$x > 2\binom{i(X)-2}{r-2} > (1 + \varepsilon)\binom{i(X)}{r-2}$$ for all sufficiently large $n$. Let us now apply Theorem \ref{cntE}. We have $$e(A) > \frac{c_0x^2}{i(X)}$$ for all sufficiently large $n$. Now we are ready to estimate the probability.

\begin{multline*}
    \sum_{\Tilde{c}n^{\frac{r-3}{r-1}} < i < \varepsilon'n}\sum_{ cin^{r-3} < x < \min(\varepsilon n^{r-2}, \binom{i}{r})} \binom{n}{2}\binom{\binom{n-2}{r-2}}{x-1}\binom{\binom{i}{r}}{x}2^{-\frac{c_0x^2}{i}} \leq \\ \leq \sum_{\Tilde{c}n^{\frac{r-3}{r-1}} < i < \varepsilon'n}\sum_{ cin^{r-3} < x < \min(\varepsilon n^{r-2}, \binom{i}{r})} 2^{\alpha x\ln n-\frac{c_0x^2}{i}}.
\end{multline*}
Note that the maximum of $$\alpha x \ln n - \frac{c_0x^2}{i}$$ is reached (for a fixed $i$) at $$x = \frac{1}{2c_0}\alpha i\ln n.$$ This means that for all sufficiently large $n$, the maximum of this value under the condition $x \geq cin^{r-3}$ is reached at the point $cin^{r-3}$ and is equal to $-c_0c^2in^{2r-6} + O(n^{r-2})$. Thus, the desired probability is estimated by the value
\begin{multline*}
    \sum_{\Tilde{c}n^{\frac{r-3}{r-1}} < i < \varepsilon'n}\sum_{ cin^{r-3} < x < \min(\varepsilon n^{r-2}, \binom{i}{r})} 2^{\alpha x\ln n-\frac{c_0x^2}{i}} \leq \\ \leq \sum_{\Tilde{c}n^{\frac{r-3}{r-1}} < i < \varepsilon'n}n^{r-2} 2^{-c_0c^2in^{2r-6} + O(n^{r-2})} \leq \\ \leq n^{r-1} 2^{-{c_0c^2\Tilde{c}n^{2r-6+\frac{r-3}{r-1}}} + O(n^{r-2})} = 2^{-{c_0c^2\Tilde{c}n^{2r-6+\frac{r-3}{r-1}}} + O(n^{r-2})} = o(1).
\end{multline*}

\section*{Acknowledgments}

I would like to thank A.M. Raigorodskii for the formulation of the problem and  the fruitful discussion of the results. I would also like to thank E.A. Yarovikova (Neustroeva) for finding the error in \cite{PyadBug} and a lot of useful comments that made the text of this paper much clearer.

\printbibliography

\end{document}